\documentclass[nofootinbib, showpacs,twocolumn,preprintnumbers,amsmath,amsfonts,amssymb,floatfix,aps,superscriptaddress]{revtex4}

\usepackage{graphicx}
\usepackage{epsfig}
\usepackage[hang,nooneline]{subfigure}
\usepackage[normalem]{ulem}
\usepackage{color}
\usepackage{hyperref}

\newcounter{fig}

\usepackage{scalerel}
\usepackage{stackengine}
\stackMath
\def\hatgap{-5.5pt}
\def\subdown{-3.3pt}
\newcommand\what[2][]{%
\renewcommand\stackalignment{l}%
\stackon[\hatgap]{#2}{%
\stretchto{%
    \scalerel*[\widthof{$#2$}]{\kern-.6pt\text{\textasciicircum}\kern-1.1pt}%
    {\rule[-0.8\textheight]{1ex}{\textheight}}
}{2ex}
_{\smash{\belowbaseline[\subdown]{\scriptscriptstyle#1}}}%
}}
\newcommand{\hats}[1]{\what[s]{#1}}

\def\XXint#1#2#3{{\setbox0=\hbox{$#1{#2#3}{\int}$}
     \vcenter{\hbox{$#2#3$}}\kern-.5\wd0}}

\begin{document}

\title{Revisiting Diffusion: Self-similar Solutions and
the $t^{-1/2}$ Decay in Initial and Initial-Boundary Value Problems}
\author{P.G. Kevrekidis}
\affiliation{Department of Mathematics and Statistics, University of Massachusetts
Amherst, Amherst, MA 01003-4515, USA}
\author{M.O. Williams}
\affiliation{Department of Chemical and Biological Engineering and PACM, Princeton University, Princeton, NJ 08544, USA}
\author{D. Mantzavinos}
\affiliation{Department of Mathematics and Statistics, University of Massachusetts
Amherst, Amherst, MA 01003-4515, USA}
\author{E.G. Charalampidis}
\affiliation{Department of Mathematics and Statistics, University of Massachusetts
Amherst, Amherst, MA 01003-4515, USA}
\author{M. Choi}
\affiliation{Department of Chemical and Biological Engineering and PACM, Princeton University, Princeton, NJ 08544, USA}
\author{I.G. Kevrekidis}
\affiliation{Department of Chemical and Biological Engineering
and PACM, Princeton University, Princeton, NJ 08544, USA}

\date{\today}

\begin{abstract}
  The diffusion equation is a universal and standard textbook model for partial differential equations (PDEs). In this work, we revisit its solutions,
  seeking, in particular, self-similar profiles. This problem connects to the
classical theory of special functions and, more specifically, to the Hermite  as well
as the Kummer hypergeometric functions. Reconstructing the solution of the original diffusion model from novel self-similar solutions of the associated self-similar PDE,
we infer that the $t^{-1/2}$ decay law of the diffusion
amplitude is {\it not necessary}. In particular,
it is possible to engineer setups of {\it both} the Cauchy problem and the initial-boundary value problem in which the solution decays
at a {\it different rate}. 
Nevertheless, we observe that the $t^{-1/2}$
rate corresponds to the dominant decay mode 
 among integrable
 initial data, i.e., ones corresponding to finite mass.
 Hence, unless the projection to such a mode
is eliminated,  generically this decay will be the
slowest one observed. 
In initial-boundary value problems,
an additional issue that arises is whether the boundary data are \textit{consonant} with the initial data; namely, whether the boundary data agree at all times with the solution of the Cauchy problem associated with the same initial data,  when this solution is evaluated at the boundary of the domain.   
In that case, the power law dictated by the solution of the Cauchy problem will be selected. 
On the other hand, in the non-consonant cases 
a decomposition of the problem into a self-similar and a
non-self-similar one is seen to be beneficial in obtaining
a systematic understanding of the resulting solution. 
\end{abstract}


\maketitle

\section{Introduction}

The model of diffusion is a textbook one both at the microscopic
level of Brownian motion~\cite{chandra},
as well as at the macroscopic PDE level and its mathematical
analysis~\cite{strauss,evans}, hence it needs no particular
introduction. In its one-dimensional form, to which we will restrict
our considerations herein, it reads
\begin{eqnarray}
  u_{t}=u_{xx},
  \label{deq1}
\end{eqnarray}
where $u=u(x, t)$ represents a physical (dependent) variable such as
temperature or concentration, and  the subscripts denote
partial derivatives with respect to time $t$ and space $x$.
This is a model so widely studied that
it is hard to envision any elements of novelty in its study
at present.

Nevertheless, the 1969 work of~\cite{cole}  identified
an apparently
previously unknown class of solutions of Eq.~\eqref{deq1} using the method of
similarity variables. When connecting with the special case
of self-similar solutions whose spatial dependence arises
in terms of the traditional self-similar variable $x/t^{1/2}$, the authors of~\cite{cole} derive
special solutions associated with parabolic cylinder functions
(a special case of the confluent hypergeometric series). In this
context, they remark
that in order to have a solution vanishing (``actually exponentially'',
as they point out) as the similarity variable tends to $\pm \infty$, these solutions must be characterized
by an integer index. They also remark that if the total
mass is constant, then the ``standard'' solution,
namely
\begin{eqnarray}
  u(x,t) \propto t^{-1/2} e^{-\frac{x^2}{4 t}},
  \label{deq2}
\end{eqnarray}
must be chosen.  The application of the method was subsequently extended
to boundary value problems~\cite{bluman2}, while other authors
extended it to a variety of different settings including, e.g., 
Schr{\"o}dinger-type equations~\cite{bluman3}, and the sine-Gordon~\cite{blumen4},
nonlinear diffusion~\cite{blumen5}, and nonlinear
Boltzmann equation~\cite{blumen6}.

The identification of similarity solutions is also
by now a textbook subject~\cite{baren,galaktionov}.
Nevertheless, a recent methodology, occasionally referred to
as MN-dynamics~\cite{betelu} (see also~\cite{clancy,siettos}, as
well as~\cite{siettos2}, where a general formulation thereof was
presented for nonlinear PDE problems) offers a simple and
systematic alternative to deriving such waveforms and the corresponding
scaling properties.
It is worthwhile to note in passing that there is a considerable volume of
literature in the study of such problems in dispersive equations, see,
e.g.,~\cite{sulem,fibich}.
It is this MN-dynamics methodology that we will adopt
here, in a prototypical problem such as the diffusion equation.

What we first obtain is in fact the self-similar solutions of~\cite{cole}, but now for {\it arbitrary} real values of the relevant index as opposed to just integer values of this index. From this, we infer that  scaling laws {\it different} from the standard $t^{-1/2}$ law of Eq.~\eqref{deq2}  are not only feasible, but actually entirely realizable
in both the initial value (Cauchy) problem and the 
initial-boundary value problem setting (in the case of finite domain considerations).
These decay laws can bear {\it arbitrary} negative exponent
in the (temporal) scaling of the solution amplitude.
It is
shown on the basis of a general class of initial conditions
and of the (full) solution of the MN-dynamics problem that the $t^{-1/2}$
decay is the slowest  one (for integrable initial data of finite
mass), and it is
explained under what conditions a different decay rate will be observed.

Similar considerations are relevant to examine when boundary conditions are present. Furthermore, it is now important to also discuss the role of conservation laws. In this case, we introduce the notions of \textit{compatible} and \textit{consonant} boundary conditions. The former ones are identified as boundary conditions that agree with the initial condition at the endpoints of the spatial domain and at time $t=0$ (but not necessarily at all times), while the latter ones correspond to boundary conditions that match at all times
the solution of the associated Cauchy problem,
when this solution is evaluated at the boundary of the domain of the initial-boundary value problem.

As we will see below, in the case of the Cauchy problem there exist initial data for which there are ``eigenfunctions'' that lead
to self-similar evolution with a specific decay rate. On the other hand, there are initial conditions that do not project solely on one such eigenfunction but have a finite projection on multiple modes, decaying
at different rates and resulting in non-self-similar evolution
with a single rate of decay. 

Similarly, in the case of initial-boundary value problems   there are choices of boundary conditions that are conducive to self-similar dynamics, and others that are not. The latter category includes  typical examples taught in undergraduate courses,
such as homogeneous Dirichlet, homogeneous Neumann, and constant
coefficient Robin conditions (cf., e.g., Chapter 4 in~\cite{strauss}).
These case examples, which are solvable by different forms of
Fourier series and bear the associated
exponential time dependence, are not inherently self-similar.
In the conducive category, 
there are time-dependent boundary conditions
that are consonant with self-similar evolution. 
For instance,
let us consider an exact eigenfunction of the self-similar MN diffusion
problem discussed above, and boundary conditions that are consonant with it (i.e., in agreement with the value
of the relevant eigenfunction at the boundaries of the domain over all times).
It is
clear that in this case the solution  is ``transparent'' to
the presence of the boundary conditions.
In this setting, we can still engineer solutions with various
kinds of decay rates (different than $t^{-1/2}$). 

Another
intriguing possibility  is that of boundary conditions
that are compatible with the initial condition (e.g. continuous so as to avoid Gibbs-type phenomena~\cite{strauss}), yet not consonant.
In this scenario, we will advocate the decomposition of the problem into (i) an initial-boundary value problem which is both consonant
and compatible, and which takes care of the boundary conditions, and (ii) a complementary
problem with homogeneous boundary data that does not have a self-similar solution.  In what follows, we will restrict our considerations to the above cases and will not examine  the definitively non-self-similar evolutions of general boundary conditions that are neither consonant, nor compatible with any self-similar waveform.

Our presentation is structured as follows. 
In Section~2, we present
the MN-dynamics approach for the diffusion equation and its basic conclusions.
In Section~3, we consider a number of select numerical experiments
in the initial-boundary value problem setting. 
Finally, in Section~4
we summarize our findings and provide a discussion of future perspectives.

\section{Theoretical Analysis: Self-similar Solutions}

As discussed in~\cite{betelu,clancy,siettos} (see also
the recent exposition in the Appendix of~\cite{siettos2}),
the scaling ansatz for seeking a self-similar solution of Eq.~\eqref{deq1} via the MN-dynamics approach 
is of the form
\begin{eqnarray}
  u(x,t)=A(\tau) w\left(\xi,\tau \right), \quad \xi = \frac{x}{L(\tau)}, \quad \tau = \tau(t),
  \label{deq3}
\end{eqnarray}
where $\xi$ is the similarity variable and the functions $\tau$, $A$ and $L$ are to be determined.

Direct substitution and division by $A$ yields
\begin{eqnarray}
  \left(w_{\tau} + \frac{A_{\tau}}{A} w - G \xi w_{\xi} \right)
  \tau_t= \frac{1}{L^2} w_{\xi \xi}.
  \label{deq4}
\end{eqnarray}
In the above equation, we have set 
\begin{equation}\label{G-def}
\frac{L_{\tau}}{L}\doteq G=\text{const.}>0,
\end{equation} 
assuming that in the $(\xi, \tau)$ frame the solution
has  a constant rate of expansion of its width
during its self-similar evolution. From this, it is immediate
to infer that 
\begin{equation}\label{L-form}
L(\tau)=L_0 e^{G \tau}.
\end{equation}
To obtain solutions that are steady in the self-similar frame,
the explicit time-dependence must be eliminated by necessitating that
\begin{eqnarray}
  \tau_t = \frac{1}{L^2}.
  \label{deq5}
\end{eqnarray}
Substituting $L$ from Eq.~\eqref{L-form} into Eq.~\eqref{deq5}
yields the relation between the old and new time frames as
\begin{eqnarray}
  e^{\tau} = \left[\frac{2 G}{L_0^2} \left(t - t^{\star} \right) \right]^{\frac{1}{2 G}}.
  \label{deq6}
\end{eqnarray}
As will be evident below, the  positive constant $G$ can be scaled out of the equations. Hence, upon manifesting this scaling, we will
select $G=1$ for simplicity. Moreover, evaluating 
Eq.~\eqref{deq6} at $t=0$ naturally reveals that $t^{\star}<0$, suggesting the well-known feature that diffusion processes blow up in reverse
time. Note that in the rescaled temporal variable, this time corresponds
to $\tau \rightarrow -\infty$.

A key observation that distinguishes the present problem
from other ones where the above method has been applied (such as,
e.g.,~\cite{siettos,siettos2}) is that now
the amplitude scaling is {\it not} fixed by the PDE dynamics
(self-similarity of the first kind~\cite{betelu}), but rather
it has to be obtained from the solution of an eigenvalue problem (see below;
this is the self-similarity of the second kind~\cite{betelu}).
In particular, we have the freedom to select 
\begin{equation}
\frac{A_{\tau}}{A}\doteq b = \text{const.},
\end{equation}
which in turn leads to 
\begin{equation}\label{A-form}
A(\tau)=A_0 e^{b \tau},
\end{equation} 
or, in the original frame via Eq.~\eqref{deq6},
\begin{eqnarray}\label{At-form}
  A(t)=A_0 \left[\frac{2 G}{L_0^2} \left(t - t^{\star} \right) \right]^{\frac{b}{2 G}}.
  \label{deq7}
  \end{eqnarray}
Upon the above scaling choices, the original PDE in the
renormalized frame, i.e., Eq.~\eqref{deq4}, can be expressed as
\begin{eqnarray}
  w_{\tau}=w_{\xi \xi} + G \xi w_{\xi}- b w.
  \label{deq8}
\end{eqnarray}
In this setting, the factor $G$ can be removed completely via rescaling space by $G^{1/2}$ and sending $b \rightarrow b/G$, hence it will be set to unity
hereafter.

Making the additional explicit transformation
\begin{equation}
w=e^{-\xi^2/2} W,
\end{equation} 
as well as slightly modifying the $(\xi, \tau)$ frame according to the change of variables
\begin{equation}
\xi \rightarrow \tilde{\xi}=\frac{\xi}{\sqrt{2}}, \quad \tau \rightarrow
\tilde{\tau}=\frac{\tau}{2},  
\end{equation} 
we obtain in the new $(\tilde \xi, \tilde \tau)$ frame the equation
\begin{eqnarray}
  W_{\tilde{\tau}}=W_{\tilde{\xi} \tilde{\xi}} -2 \tilde{\xi}
  W_{\tilde{\xi}} +2 \nu W,
  \label{deq9}
\end{eqnarray}
where 
\begin{equation}\label{nu-def}
\nu=-\left(b+1\right).
\end{equation}
The steady state problem of Eq.~\eqref{deq9} (i.e., the one obtained by setting the right hand side equal to zero)
is {\it precisely} the Hermite differential equation.
From this reduction, we can infer the self-similar solutions,
thereby connecting the results with the earlier work of~\cite{cole}.
If we assume, in particular, that in the renormalized
frame the motion is self-similar in the original
variables (i.e., steady in the new frame), then the resulting solutions are of the form
\begin{equation}
\hskip -2.2mm 
w=e^{-\frac{\xi^2}{2}} \!\left[ c_1(\nu) H_\nu\!\left(\frac{\xi}{\sqrt{2}}\right)
  \! + c_2(\nu) {}_1F_{1}\!\left(\!-\frac{\nu}{2},\frac{1}{2},\frac{\xi^2}{2}\right)\! \right]\!, 
  \label{deq10}
\end{equation}
where $H_\nu$ denotes the Hermite functions, while $_1F_1$ denotes the
so-called Kummer confluent hypergeometric series
functions. These solutions are also related to the
so-called Weber or parabolic cylinder functions (see~\cite{ww} for more details). 

It is now important to clarify some points regarding the
solutions \eqref{deq10}. Firstly, the Kummer functions are {\it even} for any value of $\nu$. On the other hand, the Hermite functions have definite parity only in the case of \textit{integer} $\nu$. Actually, in that case they reduce to the well-known Hermite polynomials, which are even for $\nu$  even and odd for $\nu$  odd.
Furthermore, importantly, and {\it differently} from what was suggested
in~\cite{cole}, $\nu$ does {\it not} have to be an integer.
Indeed, the solutions of Eq.~\eqref{deq8} are still well-defined  functions for any $\nu$ (equivalently, $b$)  real. However, as may be evident
when interpolating between, e.g., an even and an odd value of $\nu$, the Hermite
solution $H_\nu$ for general $\nu$ is an \textit{asymmetric} one. 

Another
key point to make is that bounded solutions exist only for
$b<0$ (equivalently, for $\nu>-1$), while integrable solutions exist only for
$b<-1$ (equivalently, for $\nu>0$). These distinctions will be important in what follows, not only for mathematical reasons but also for physical ones,
since integrability here is tantamount to finite mass, a requirement
of particular physical relevance.
Yet another important observation concerns asymptotic properties. In the
case of integer $\nu$, the standard asymptotic properties
that we expect from the Hermite polynomials in connection to
$w$ apply i.e., asymptotically the solutions decay
as $w \sim \xi^\nu e^{-\xi^2/2}$; this is the case referred
to in~\cite{cole}. {\it However}, in the case where this
integer ``quantization'' of $\nu$ is absent, the stationary series solution of Eq.~\eqref{deq9} does not close as a regular, finite-order polynomial. Instead, it produces a more general function and, in this
case, it is known~\cite{lebedev} that, e.g., for the Kummer
solution the asymptotics as the argument tends to infinity
are 
\begin{eqnarray}
  {}_1F_{1}(\alpha, \beta; z)=e^z z^{\alpha-\beta} \frac{\Gamma(\beta)}{\Gamma(\alpha)}.
  \label{deq11}
\end{eqnarray}
Straightforward substitution reveals that in this
general, non-quantized case the asymptotics decay in the
form of a power law according to $w \sim \xi^{b}$,
once again suggesting that for decay, we need $b<0$, while for integrability we must have $b<-1$. 

Another important observation concerns the cases where $\nu$ is an even integer. In this
case, the two solutions of Eq.~\eqref{deq8} {\it coincide}.
This is natural to expect since for these particular
values the Hermite polynomials are even and, as mentioned
previously, the Kummer functions are always even. Hence, a second
linearly independent solution is present for $\nu$ even.
This can be found on a case-by-case basis.
For example, for $\nu=0$ ($b=-1$) this is $w=e^{-\xi^2/2} {\rm Erfi}(\xi/\sqrt{2})$
(where {\rm Erfi} is the imaginary error function).
Similarly, for $\nu=2$ ($b=-3$) we have 
$w = 2 \xi + \sqrt{2 \pi} \left(1-\xi^2\right) e^{-\xi^2/2} {\rm Erfi}(\xi/\sqrt{2})$, and so on.

\onecolumngrid

\begin{figure}[tbp]
\begin{center}
\mbox{\hspace{-0.5cm}
\subfigure[][]{ 
\includegraphics[height=.24\textheight, angle =0]{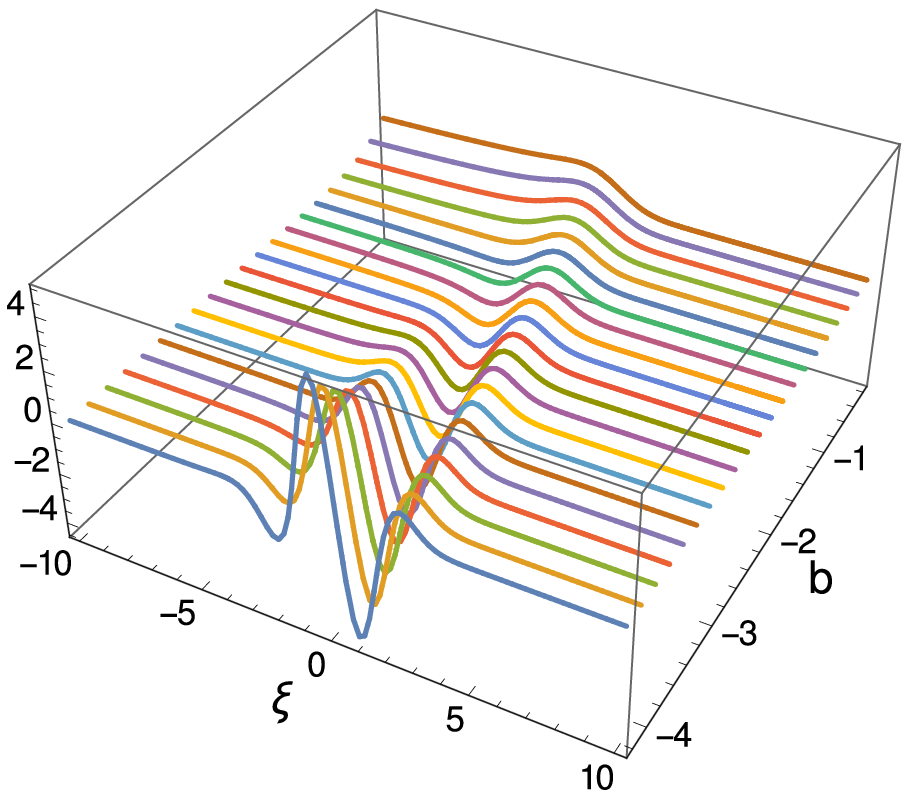}
\label{diff1_c}
}
\subfigure[][]{\hspace{1.5cm}
\includegraphics[height=.24\textheight, angle =0]{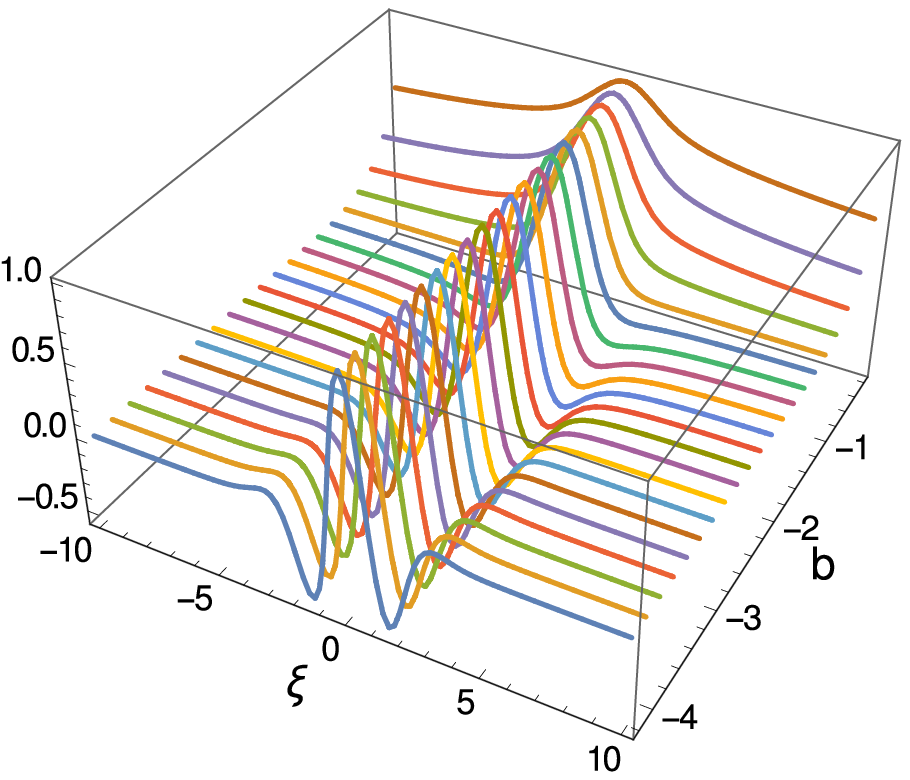}
\label{diff1_d}
}
}
\mbox{\hspace{-0.0cm}
\subfigure[][]{ 
\includegraphics[height=.26\textheight, angle =0]{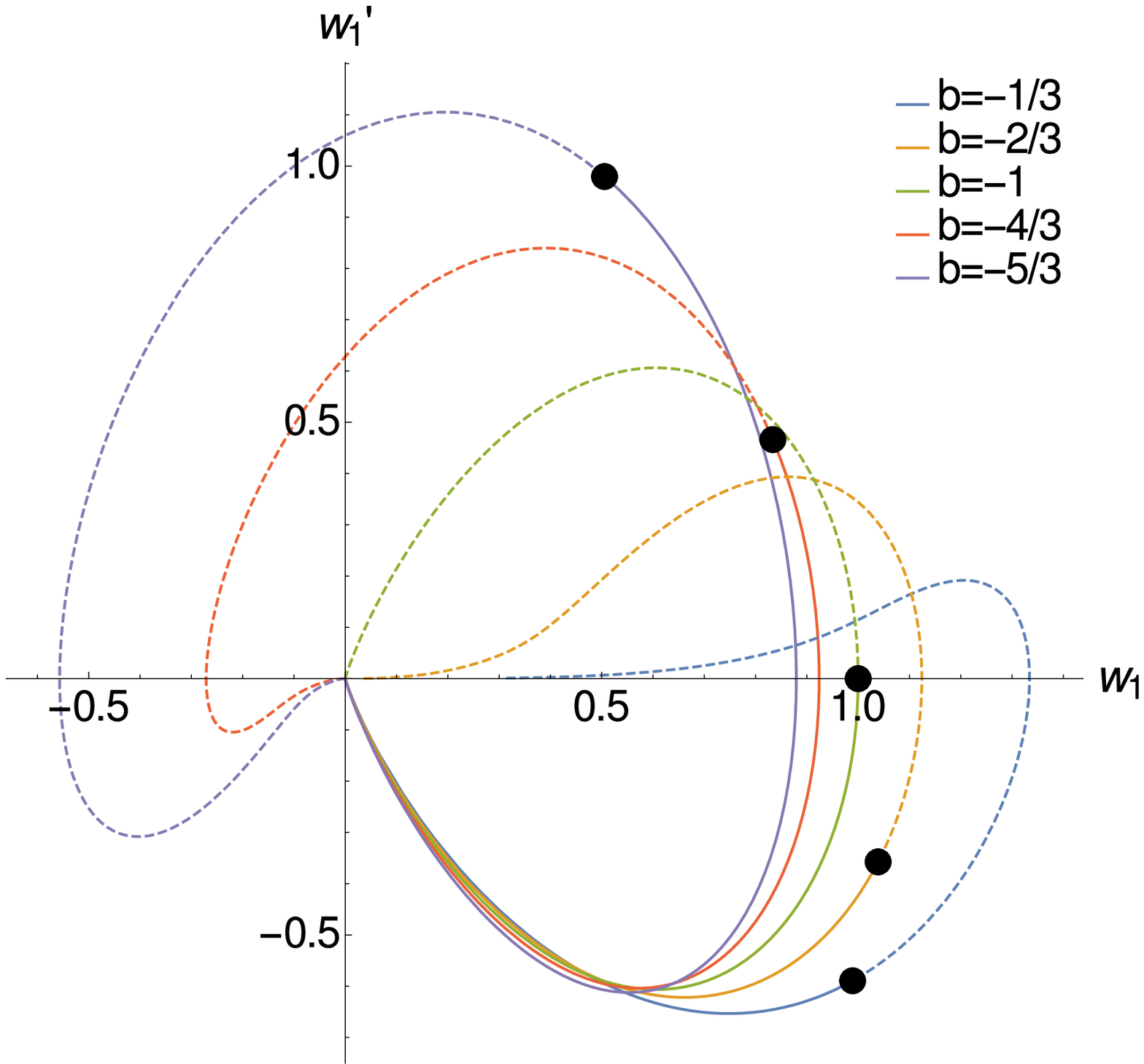}
\label{diff1_e}
}
\subfigure[][]{\hspace{1cm}
\includegraphics[height=.26\textheight, angle =0]{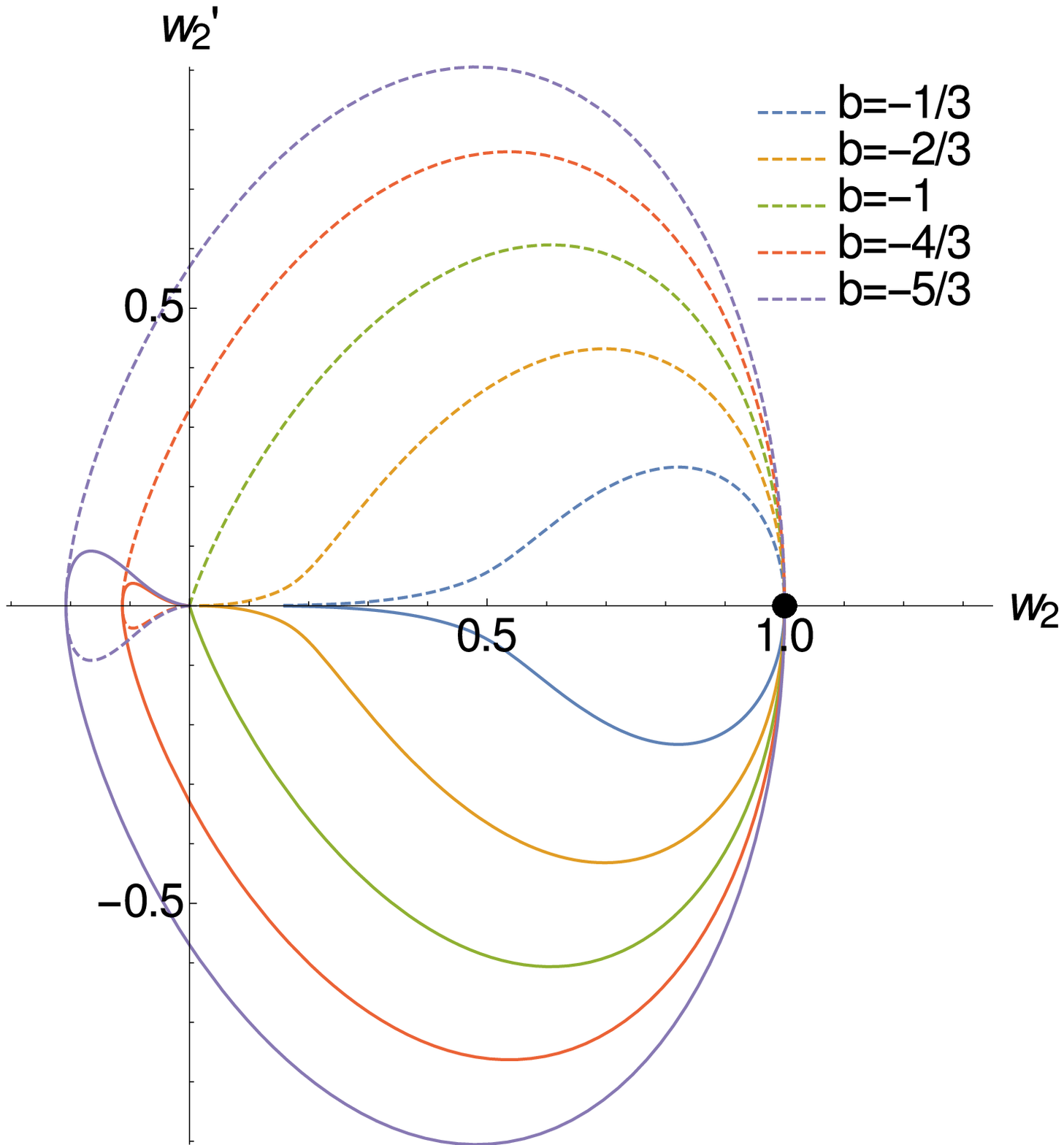}
\label{diff1_f}
}
}
\end{center}
\vskip -5mm
\caption{\textit{Left} and \textit{right} columns respectively correspond to the Hermite component $w_1 = e^{-\xi^2/2} H_{\nu}\left(\xi/\sqrt 2\right)$ and the Kummer component $w_2 = e^{-\xi^2/2}{}_1F_{1}\left(-\nu/2, 1/2, \xi^2/2\right)$
(recall $\nu = -b-1$) of the steady state, self-similar solution~\eqref{deq10}. 
(a) and (b): Plots of $w_1$ and $w_2$ against both $\xi$ and $b$. The thicker curves correspond to integer values of $b$.
(c) and (d): Phase portraits, i.e. plots of $w_{1,2}'$ against $w_{1,2}$, 
for $\xi \in (-\infty,0] \text{(dashed)} \cup [0,\infty)
\text{(solid)}$ and various values of $b$. 
The black dots correspond to $\xi=0$.
}
\label{diff1}
\end{figure}

\twocolumngrid

The general features of the Hermite and Kummer components of the steady state, self-similar solution~\eqref{deq10} of Eq.~\eqref{deq8} discussed above are visualized in Fig.~\ref{diff1}. 
In the left column of  Fig.~\ref{diff1}, we provide plots of the Hermite component against $\xi$   for various choices of  $b=-\left(\nu+1\right)$.
The same is done for the Kummer component in the right column of Fig.~\ref{diff1}. Note that for
the cases shown with $b> -1$, the non-integrability
of the wavefunction is a result of the slow decay in both figures.
Moreover, note that for $b<-1$ the solution also acquires
negative values, a feature somewhat atypical for diffusional
dynamics when considering the model, e.g., for chemical concentrations.

Armed with the above special (self-similar) solutions, we can
actually solve the original problem of Eq.~\eqref{deq1}
in the basis most naturally tailored to address self-similar
evolution profiles, namely the basis of stationary solutions~\eqref{deq10} of
Eq.~\eqref{deq8}. To do so in a more general form,
we separate variables
in Eq.~\eqref{deq8} using $w(\xi,\tau)=X(\xi) T(\tau)$.
The temporal part  yields directly an exponential decay of the form 
\begin{equation}
T(\tau)\sim e^{-\lambda \tau},
\end{equation}
where $\lambda$ is a suitable eigenvalue, while the (rescaled) spatial part  $X(\xi)$ satisfies the
same ODE as the steady state problem of \eqref{deq8}, but
now with the substitution 
\begin{equation}\label{btil-def}
b \rightarrow \tilde{b}=b-\lambda.
\end{equation}
Hence, using the separated variables solutions we can \textit{formally}
construct a linear superposition of the solutions
of Eq.~\eqref{deq8} as
\begin{align}
&w(\xi,\tau)=e^{-\frac{\xi^2}{2}}   \int_{\lambda}  e^{-\lambda \tau} 
\nonumber\\
&\times \left[ c_1(\lambda) H_{\tilde \nu}\left(\frac{\xi}{\sqrt{2}}\right)
  + c_2(\lambda)\, {}_1F_{1}\left(-\frac{\tilde \nu}{2},\frac{1}{2},\frac{\xi^2}{2}\right) \right] d\lambda,
  \label{deq12}
\end{align}
where  
\begin{equation}\label{nutil-def}
\tilde{\nu}=-\left(b-\lambda+1\right).
\end{equation}

The superposition \eqref{deq12} is done at a {formal} level, since the range of the relevant parameter $\lambda$  and the associated decay and smoothness properties of $c_1(\lambda), c_2(\lambda)$ have not been specified. 
In fact,  for solutions lying in most function spaces it is unclear what the range of integration in \eqref{deq12} or the properties of the kernel functions $c_1(\lambda), c_2(\lambda)$ should be.
In the case of e.g. $\lambda$ spanning over all real numbers, it is not known whether the ``component'' functions $H_{\tilde \nu}(\xi/\sqrt 2)$ and ${}_1F_1(-\tilde \nu/2, 1/2, \xi^2/2)$ form a well-defined basis  (with an appropriate inner product etc.). 

On the other hand,  it can be shown~\cite{nw09} (see also~\cite{gw2012})
that   for exponentially decaying initial data, there exists a Hilbert space such that the superposition Eq.~\eqref{deq12} ``collapses'' in such a way that only 
the ``quantized'' modes associated with $\tilde{\nu} \in\mathbb Z$ survive, and that the resulting sum gives a convergent representation of the solution for all times.
In that case,  spectral projections and decompositions are well-defined and the solution of Eq.~\eqref{deq8} is indeed given by a variant
of  Eq.~\eqref{deq12} expressed as a sum over these integer indices. 
From a function-analytic perspective, an understanding of how to
incorporate the continuum of power law decay solutions into the
superposition of the quantized, integer-indexed, exponentially decaying ones
represents a particularly timely and relevant direction for future work.
In our numerical computations in the following section, we will
focus chiefly on rapidly decaying initial data. Nevertheless, we note
in passing that we have confirmed
that, for suitable choices of initial conditions
and corresponding consonant boundary conditions, power law solutions of the type identified above can be observed
for all time.

We now return to the  frame of the original variables $(x,t)$
in order to reconstruct in that frame the superposition appearing in Eq.~\eqref{deq12}.
Combining Eqs.~\eqref{L-form},~\eqref{deq6},~\eqref{A-form} and~\eqref{deq12} with Eq.~\eqref{deq3}, we thereby arrive at the solution formula
\begin{align}
u(x,t) &= e^{-\frac{x^2}{4 (t - t^{\star})}}
    \int_{\lambda} \left(t-t^{\star}\right)^{\frac{b-\lambda}{2}}
    \bigg[c_1 (\lambda) H_{\tilde{\nu}}\left(\frac{x}{2 \sqrt{t-t^{\star}}}\right)
\nonumber\\
&\quad
+ c_2(\lambda)\, {}_1F_{1}\left(-\frac{\tilde{\nu}}{2},\frac{1}{2},\frac{x^2}{4
    \left(t-t^{\star}\right)}\right) \bigg] d\lambda.
\label{deq13}
\end{align}

The superposition~\eqref{deq13} sheds some direct light to
the selection of the $t^{-1/2}$ amplitude decay. In particular, recall (see earlier discussion about the behavior of the stationary solution \eqref{deq10}) that in order for $w$ as given by Eq.~\eqref{deq12} to manifest decay as $|\xi|\to \infty$, it must be  that $\tilde b<0$,
while $\tilde b<-1$ must be satisfied for $w$ to be integrable with respect to $\xi$.
Thus,  among all integrable solutions of the form \eqref{deq13} (i.e., from a physical perspective,
among all solutions of finite mass), the limiting case 
$\tilde b=-1^- \Leftrightarrow \tilde{\nu}=0^-$ 
corresponds to the dominant in time ``mode'', since it is the one that has the slowest decay of $t^{-1/2}$. 
It is therefore this slowest mode that we should generically expect to observe in physical measurements. This is analogous, in a way, to the textbook
case that we are well familiar with from separation of variables
of diffusion with homogeneous  boundary conditions
on a finite interval. For example, in the case of homogeneous Dirichlet conditions on the interval $[0, D]$, the solution reads
\begin{eqnarray}
  u(x,t)=\sum_{n \in {\mathbb Z}} A_n\, e^{- \left( \frac{n \pi}{D}\right)^2 t}
    \sin\left(\frac{n \pi x}{D}\right),
    \label{deq14}
\end{eqnarray}
and for generic initial data bearing a projection to the dominant
($n=1$) mode of Eq.~\eqref{deq14}, all additional terms
die out (in this latter case
exponentially fast) and only the ground state proportional
to $\sin(\pi x/D)$ survives, decaying at the rate of the
slowest mode.

However, although our analysis suggests that, generically, in the
Cauchy problem the dominant  mode will decay like  $t^{-1/2}$,
at the same time we explicitly illustrate that it is possible to
prescribe {\it any} rate of decay at will, provided that a suitable initial condition is prescribed. For instance,
in the case of an initial condition for which the coefficient $c_1(\lambda)$ is a Dirac $\delta$-function 
centered at $\lambda=b+2$  (thus implying $\tilde{\nu}=1$ through Eq.~\eqref{nutil-def}),
only the mode $H_1$ arises, whose decay rate is {\it different}
from $t^{-1/2}$. In particular, for $u(x,0)=x e^{-x^2/2} = f(x)$,
we obtain the solution
\begin{eqnarray}
  u(x,t)=\frac{x}{(2 t+1)^{\frac 32}} e^{-\frac{x^2}{2(2 t+1)}} \equiv
  \frac{1}{2 t+1} f(\xi),
  \label{deq15}
\end{eqnarray}
where $\xi = x/\sqrt{2t+1}$, which is full agreement with the
prediction of the MN-dynamics. Here, the vanishing projection
of the initial condition on the lowest order mode enables
the observation of a different rate of decay. One may wonder
how this possibility is compatible with the well-known fact \cite{evans} that the unique solution of the Cauchy problem of Eq.~\eqref{deq1} is  given by the formula 
\begin{eqnarray}
  u(x,t)=\frac{1}{\sqrt{4 \pi t}} \int_{y\in\mathbb R}
    e^{-\frac{(x-y)^2}{4 t}} f(y) dy.
    \label{deq16}
\end{eqnarray}
Actually, Eq.~\eqref{deq16}  produces precisely
the result of Eq.~\eqref{deq15}
upon its application to the corresponding initial condition.
Hence, although the fundamental solution of diffusion may decay at the $t^{-1/2}$ rate, this is by no means necessary for all self-similarly decaying solutions of the Cauchy problem (and for the analogous solutions of initial-boundary value problems, as we will see in the following section).

So far, we have inferred that   the diffusion equation, which is a self-similar PDE, can manifest decay at different power law rates when supplemented with suitable ``self-similar'' initial data (i.e., data giving rise to self-similar solutions, along the directions of the ``eigenfunctions'' of the linear operator problem  considered
above). This is indeed
inferred by the explicit solution of the
Cauchy problem~\eqref{deq16}. We now turn to the examination of how
this conclusion is going to be affected by the presence of
boundary conditions.

\section{Numerical Results: Initial-Boundary Value Problems}

\begin{figure}[tbp]
\centering
\hskip 1.3cm
\includegraphics[width=6cm]{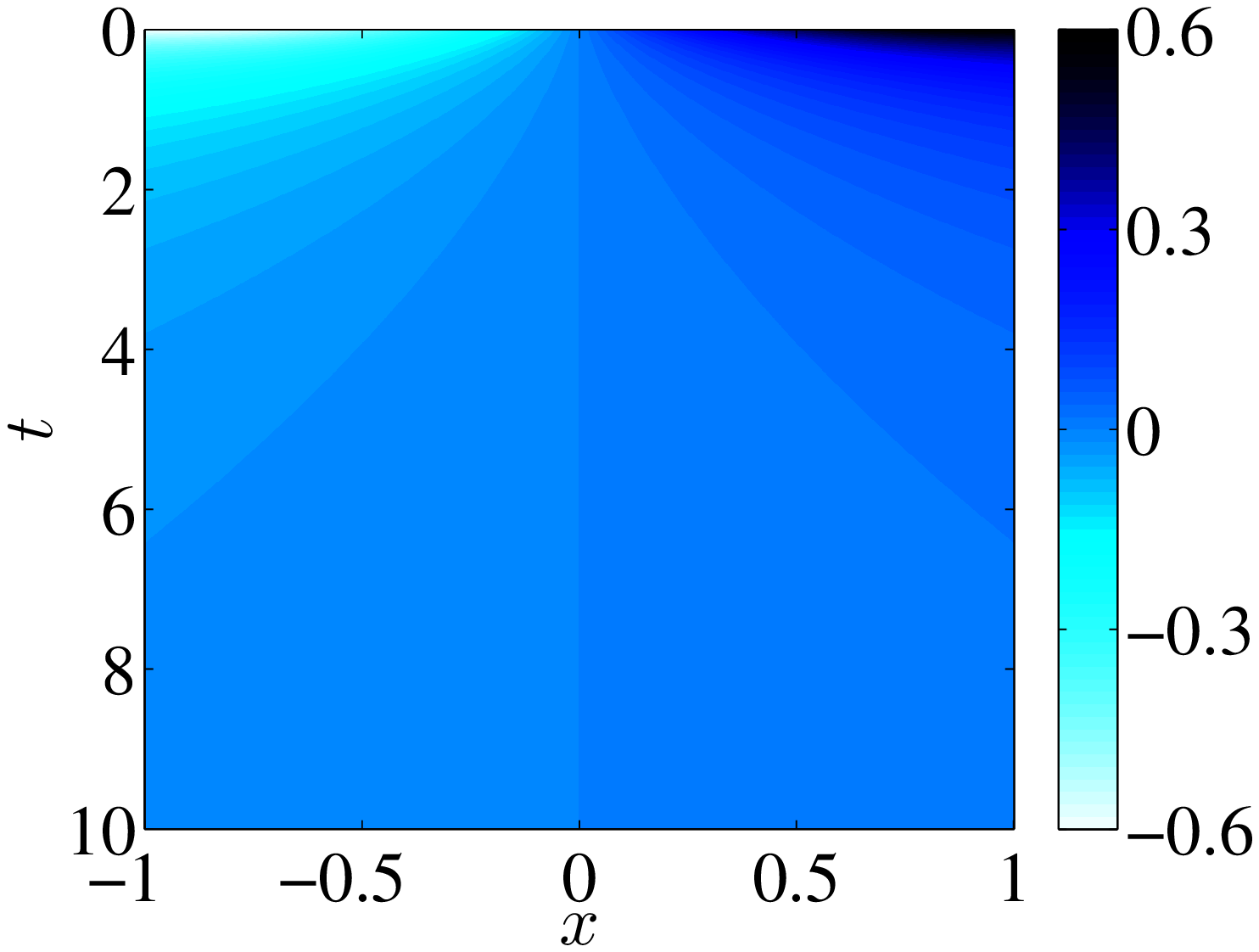}
\includegraphics[width=6cm]{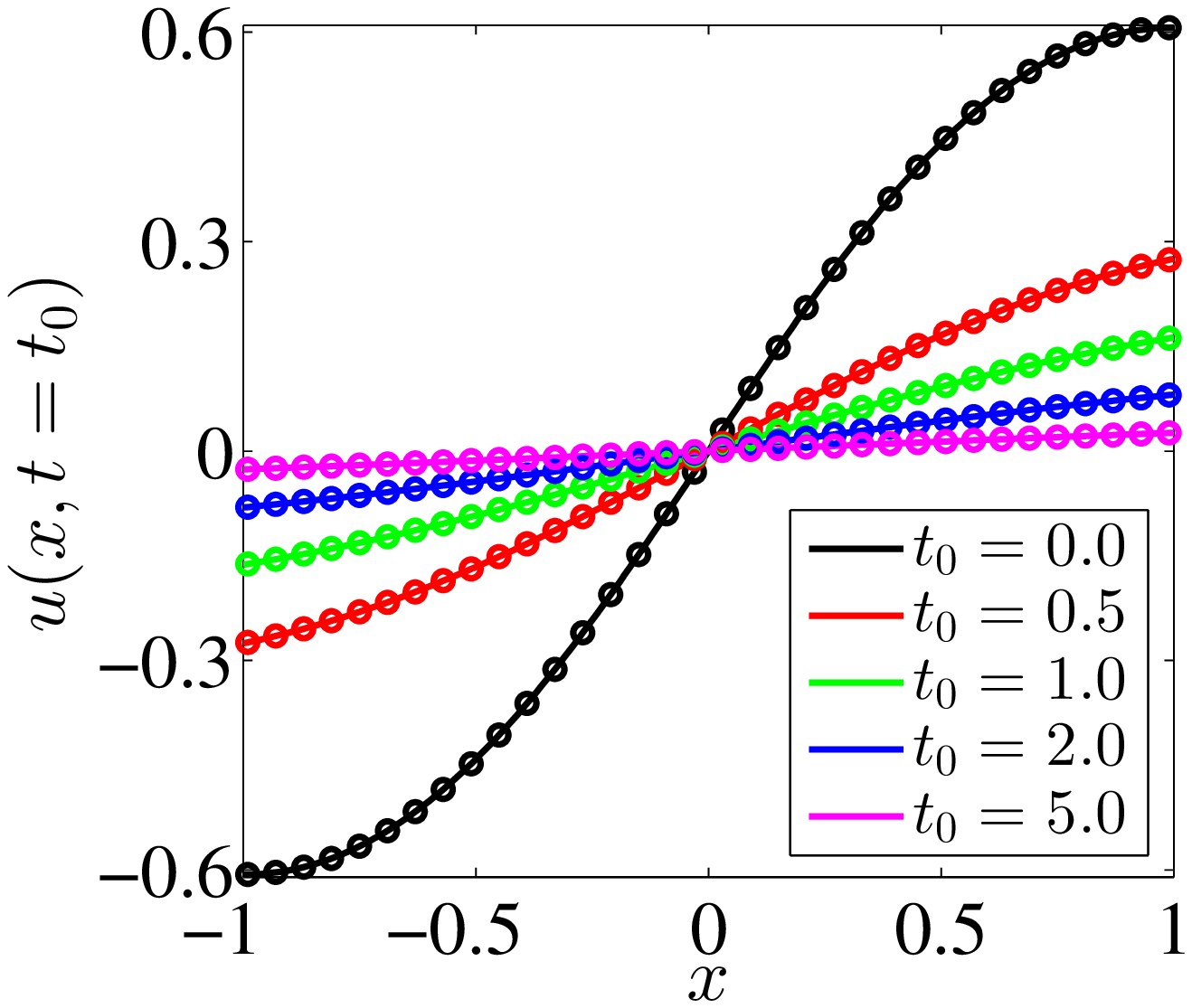}
\includegraphics[width=6cm]{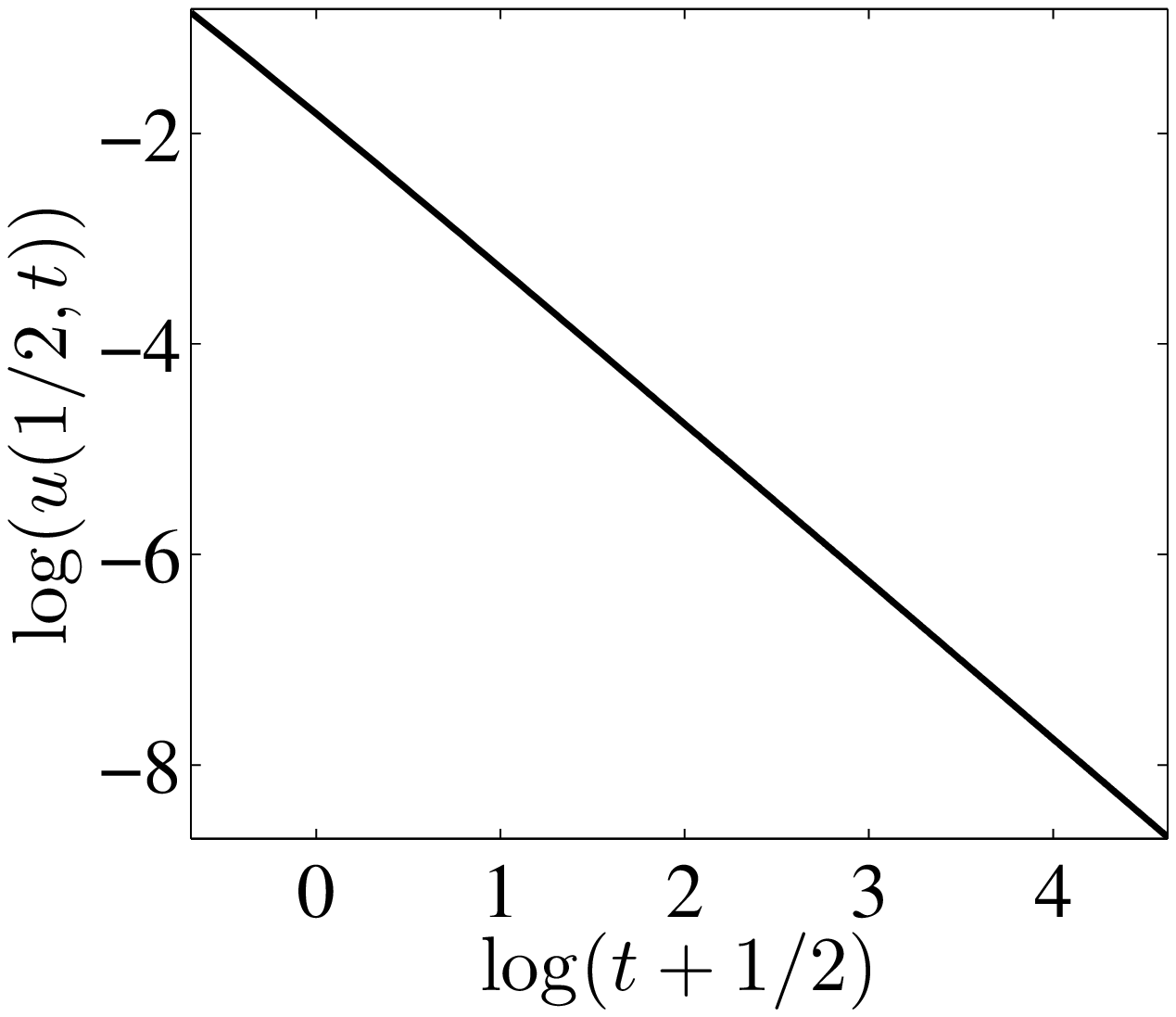}
\caption{
Spatiotemporal evolution (top panel) and spatial profiles 
at different times $t=t_0$ (middle panel) of the solution of the diffusion 
equation~\eqref{deq1} on the interval $[-D, D]$ with initial condition $u_0(x) = x e^{-x^2/2}$, 
boundary conditions $u(\pm D, t) = \pm \frac{D}{(2t+1)^{3/2}} e^{-\frac{D^2}{2(2t+1)}}$
and $D=1$. The boundary conditions are both compatible and consonant with 
the initial condition; hence, the solution is simply the self-similar solution
associated with the initial condition in the Cauchy problem setting
(cf. Eq.~\eqref{deq16}, and 
its (asymptotic for large time)
decay rate is $t^{-3/2}$ as is shown in the bottom panel of the figure.
} \label{case4}
\end{figure}

We now turn our attention to the case where boundary conditions
are also present. Continuing upon the discussion at the end of the
previous section, properties of Eq.~\eqref{deq1} such as the maximum
principle, the comparison principle etc. naturally lead to the presence
of a unique solution under standard (e.g. Dirichlet, Neumann, etc.)
boundary conditions. This immediately suggests
how the decay rates (distinct from $t^{-1/2}$) identified in the
previous section can be observed in  initial-boundary value problems (IBVPs). In particular, suppose that we prescribe
a self-similar solution among those determined in the previous section
and supplement it with boundary conditions that are consonant with this
solution. For instance, continuing our example from the
previous section, let us prescribe the initial condition
$u(x,0)=x e^{-x^2/2}$ and boundary conditions  in the symmetric domain
$[-D,D]$ of the form $u(\pm D,t)= \pm \frac{D}{(2 t+1)^{3/2}} e^{-\frac{D^2}{2(2 t+1)}}$, which have been chosen to match the solution~\eqref{deq15} of the  Cauchy problem associated with the above initial condition.\footnote{Note that in this notion of consonant boundary conditions, one does not  necessarily need to prescribe Dirichlet boundary conditions that agree
  with the solution; other types, including (inhomogeneous)
  Neumann or Robin conditions, can be engineered in a similarly straightforward manner.}
Then, the unique solution of this IBVP (i.e., of Eq.~\eqref{deq1} supplemented with the above initial and boundary conditions) will obviously once again be given by  Eq.~\eqref{deq15} 
(see  Fig.~\ref{case4}). Similarly to this case example of exponentially decaying self-similar profiles, under suitable initial and boundary conditions we have also confirmed  (data not shown) the existence of power law solutions of Eq.~\eqref{deq9} for non-integer values of $\nu$.

That is to say, we can {\it engineer} variants of the IBVP
for which any rate of decay is possible.
In these cases, we can think of the boundary as being
effectively ``transparent'' to the profile of the solution,
enabling the structure to maintain the decay rate that would
be observed in the absence of the boundary i.e., in the
Cauchy problem, in the same vein as discussed above.

A question related to the presence of the boundaries, also
touched upon in~\cite{cole}, is that of conservation laws
(and, more generally, of moment equations for moments of the PDE solution). 
Arguably, this is especially important for a model such as diffusion, for which  the mass  $M=\int_{-D}^D u dx$ obeys the conservation law
\begin{eqnarray}
  \frac{dM}{dt}= u_x(D,t)-u_x(-D,t).
  \label{deq17}
\end{eqnarray}
Bearing in mind that the flux in this case is proportional
to the (opposite of the) gradient, this suggests that the
rate at which mass changes in the domain depends on how
much influx (or outflux) occurs from the boundaries (at $x=\pm D$).
Here, again, a nontrivial difference with the work of~\cite{cole}
arises. In particular, for an anti-symmetric solution such
as that of Eq.~\eqref{deq15}, the derivative $u_x$ will be
symmetric, hence $dM/dt=0$ by construction. The mass will
thus be {\it conserved} if supplemented with consonant  boundary data of the Neumann type i.e.,
data that are obtained from the derivative of the solution
evaluated at $x=\pm D$. Nevertheless, the solution will
not decay with $\tilde{\nu}=0$ (cf. Eq.~\eqref{nutil-def}), contrary to what is suggested
in~\cite{cole}.

Hence, two important remarks so far are that (i) arbitrary
decay rates of the solution {\it can} be prescribed even
for IBVPs, yet (ii) they do not necessarily conflict
with the conservation of mass. It may well be that there is
a mass inflow and outflow and yet the balance thereof enables
the conservation of the mass within the bounded domain.
Nevertheless, it is clear that the MN-dynamics solutions
and initial conditions consonant with them will necessitate
time-dependent boundary conditions and hence require a
time-dependent flux of mass through each one of the domain boundaries
(even if this does not lead to a net flux).

In light of these remarks, Neumann conditions will control
via Eq.~\eqref{deq17} the mass flow in the system. 
Nevertheless, any type of boundary condition can be made to be consonant with self-similar evolution. This includes Dirichlet 
boundary conditions, enforcing e.g. $u(\pm D,t)=A(t) w(\pm D/L(t))$,
or Robin boundary conditions,  in which case
consonance (with an exact solution $w$) requires
\begin{eqnarray}
  u_x(\pm D, t)= \frac{w'\left(\dfrac{\pm D}{\sqrt{2 \left(t-t^{\star}\right)}}\right)}{w\left(\dfrac{\pm D}{\sqrt{2 \left(t-t^{\star}\right)}}\right)} \, \frac{u(\pm D,t)}{\sqrt{2 \left(t-t^{\star}\right)}}.
  \label{deq18}
\end{eqnarray}

However, we now turn to a more intriguing scenario. In particular,
while we can envision initial and boundary conditions consonant with the
standard type of decay (i.e., $t^{-1/2}$), as well as consonant with non-standard
types of decay (i.e., with a different exponent), it is also possible
to consider non-consonant scenarios. In particular, it
is possible to initialize Gaussian data along with
boundary conditions that are compatible with it (i.e.,
time-dependent in a way that is continuous between initial
and boundary conditions at $t=0$), yet
inducing flux at a rate  associated with a different
power law decay. It is also possible to initialize, e.g., in accordance
with the solution of Eq.~\eqref{deq15},  yet use boundary data
that are consonant with $t^{-1/2}$ decay. It is then a natural question
to inquire what happens in such cases, i.e., which rate of
self-similar decay is observed (if any). It should be noted
again here that we choose the initial data to be compatible
to avoid pathologies such as Gibbs-type phenomena.

We thus now examine this scenario of boundary conditions that are
compatible but non-consonant with the initial condition. More
specifically, let us consider Eq.~\eqref{deq1} with initial condition $u(x, 0) = e^{-x^2/2}$, 
which is associated with the standard decay rate of $t^{-1/2}$, and boundary conditions  
\begin{equation}
u(\pm D, t) = \frac{c^\star}{t+1} \ {}_1F_1\left(-\tfrac 12, \tfrac 12, \tfrac{D^2}{4(t+1)}\right) e^{-\frac{D^2}{4(t+1)}},
\label{kaz1}
\end{equation}
which is associated with a self-similar decay at a rate of $t^{-1}$. Importantly, the value of the constant $c^\star$ is chosen to enforce compatibility at $x=\pm D$ and  $t=0$, namely we have
\begin{equation}
e^{-\frac{D^2}{2}}
=
c^\star  {}_1F_1\big(\!-\tfrac 12, \tfrac 12, \tfrac{D^2}{4}\big) e^{-\frac{D^2}{4}}.
\end{equation}

Our approach in solving this broad class of compatible yet non-consonant
problems is to exploit linearity in order to decompose them 
into two sub-problems: one which is {\it both} consonant and compatible, and thus exhibits the self-similar decay imposed
by the boundary conditions,  and one with homogeneous boundary data
(and appropriately modified initial data) which do {\it not} feature
self-similar decay. In this way, we suggest that the solution at all
times maintains a self-similar and a non-self-similar part, with the
latter being described by an appropriate Fourier series solution.

More specifically, we begin by writing
\begin{equation}
e^{-\frac{x^2}{2}} = u_0^{(1)}(x) + u_0^{(2)}(x),
\end{equation}
where
\begin{align}
u_0^{(1)}(x) &\doteq e^{-x^2/2} - u_0^{(2)}(x),
\\
u_0^{(2)}(x) &\doteq c^\star {}_1F_1\left(-\tfrac 12, \tfrac 12, \tfrac{x^2}{4}\right) e^{-\frac{x^2}{4}}.
\end{align}
Exploiting linearity, we then have  
\begin{equation}
u(x, t) = u^{(1)}(x, t) + u^{(2)}(x, t),
\label{kaz2}
\end{equation}
where $u^{(1)}(x, t)$ satisfies the \textit{homogeneous Dirichlet} IBVP
\begin{subequations}\label{homo-ibvp}
\begin{align}
&u^{(1)}_t = u^{(1)}_{xx}, \hskip 1.62cm x\in (-D, D), \ t>0,\\
&u^{(1)}(x, 0) = u_0^{(1)}(x), \quad x\in [-D, D],\\
&u^{(1)}(\pm D, t) = 0, \hskip .98cm t\in [0, \infty),
\end{align}
\end{subequations}
and $u^{(2)}(x, t)$ is the solution of the (self-similar, i.e.,
both consonant and compatible) \textit{non-homogeneous Dirichlet} IBVP
\begin{subequations}\label{kummer-ibvp}
\begin{align}
&u^{(2)}_t = u^{(2)}_{xx}, \hskip 2.16cm x\in (-D, D), \ t>0,\\
&u^{(2)}(x, 0) = u_0^{(2)}(x),  \hskip 0.9cm x\in [-D, D],\\
&u^{(2)}(\pm D, t) = u(\pm D, t), \quad t\in [0, \infty).
\end{align}
\end{subequations}

\begin{figure}[tbp]
\begin{center}
\includegraphics[width=6.5cm]{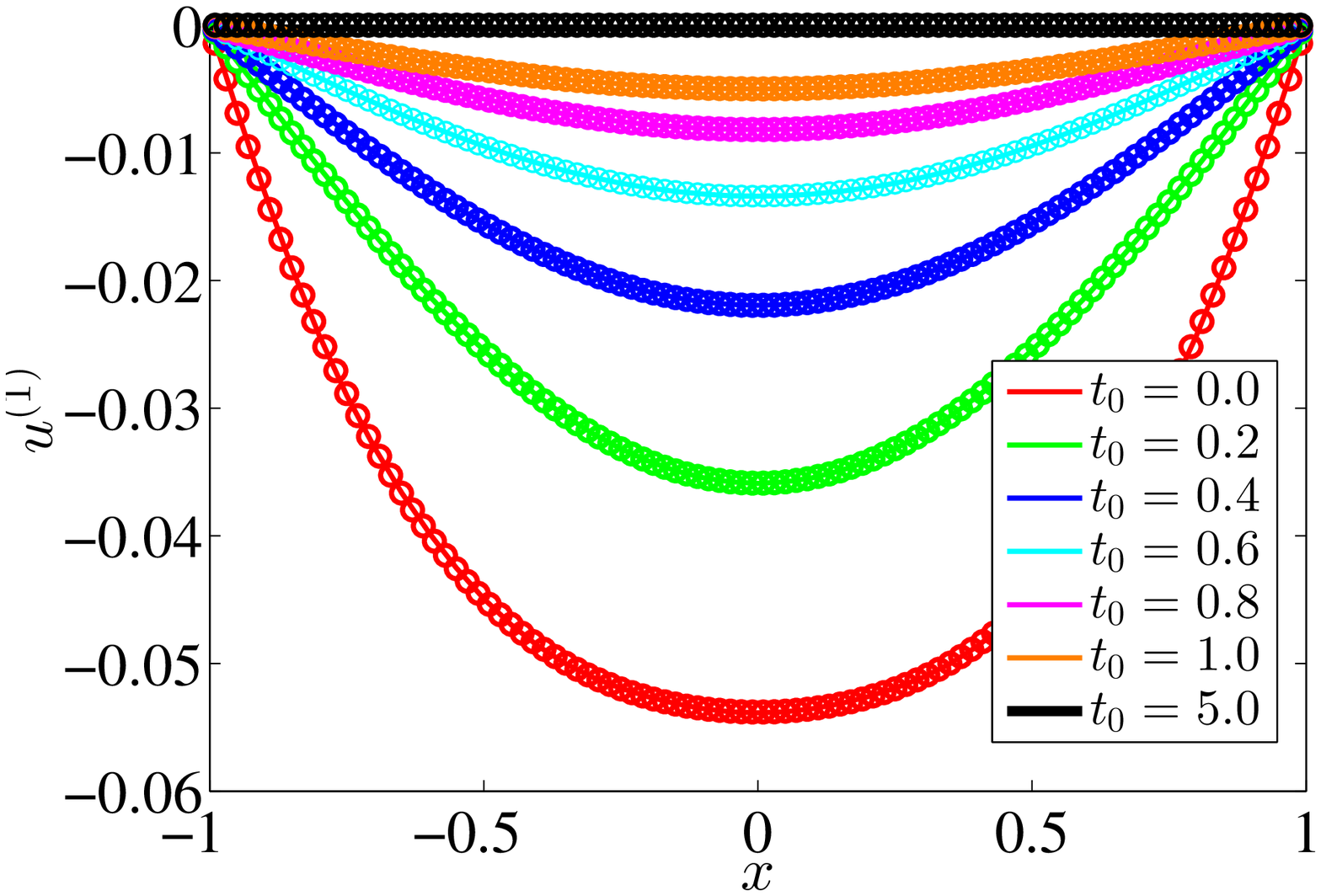}
\includegraphics[width=6.5cm]{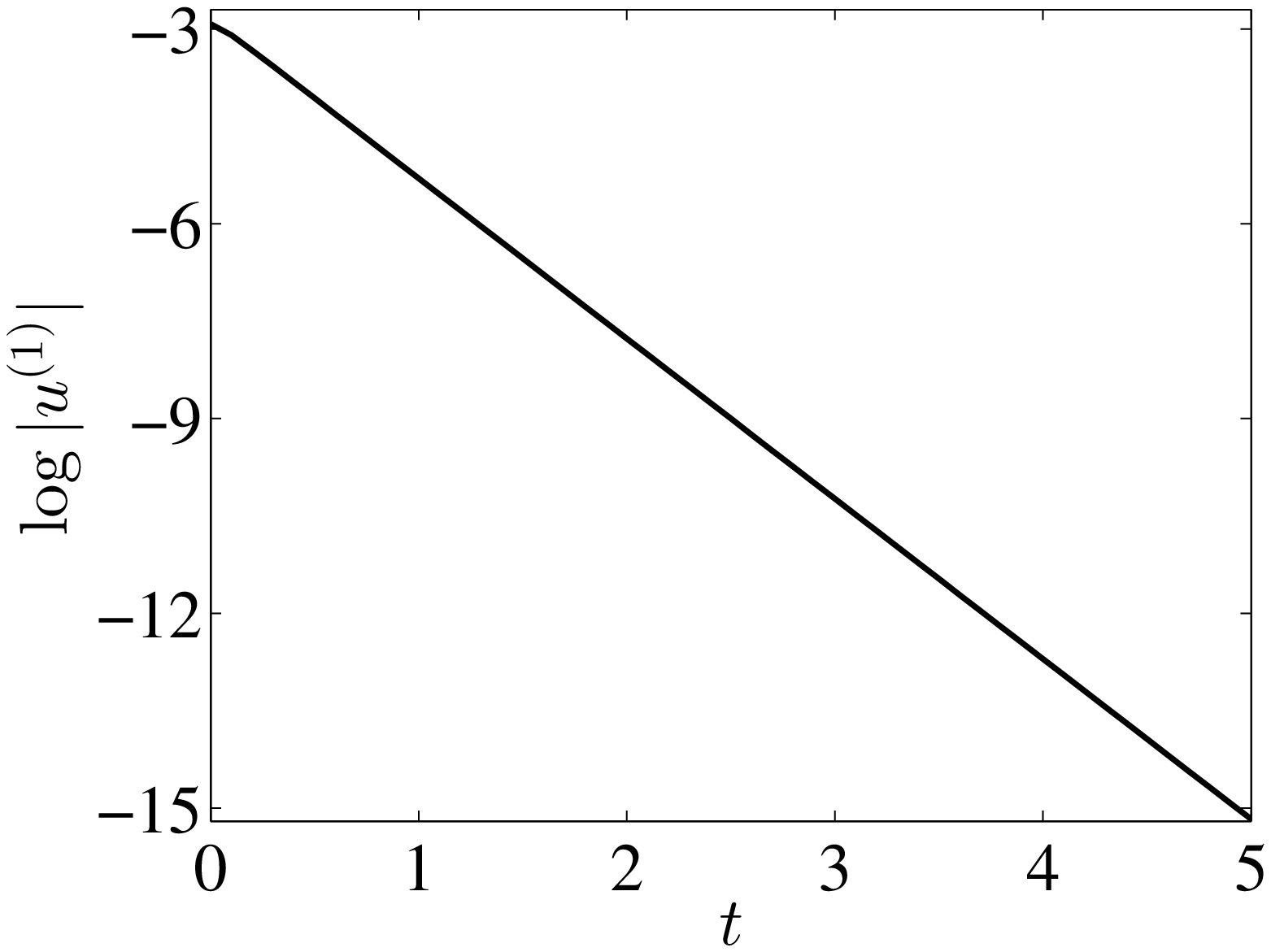}
\caption{
\textit{Top panel}: Spatial profiles of the solution $u^{(1)}$ of the homogeneous Dirichlet IBVP~\eqref{homo-ibvp} 
at various times $t=t_{0}$, obtained in two ways: 
(i) by using the first ten modes of the series representation~\eqref{utm-series} (solid lines); 
and (ii) by directly solving the IBVP numerically (open circles). 
\textit{Bottom panel}: Plot of $\log|u^{(1)}(0, t)|$. As expected, the solution decays exponentially to zero with $t$.}
\label{fig2}
\end{center}
\end{figure}

The choice of initial and boundary conditions 
in IBVP~\eqref{kummer-ibvp} implies that the solution of this problem,
due to uniqueness, must be the self-similar solution with decay rate  $t^{-1}$, i.e.
\begin{equation}\label{kaz3}
u^{(2)}(x, t) = \frac{c^\star}{t+1} \ {}_1F_1\left(-\tfrac 12, \tfrac 12, \tfrac{x^2}{4(t+1)}\right) e^{-\frac{x^2}{4(t+1)}}.
\end{equation}
On the other hand, in the case of the homogeneous Dirichlet IBVP~\eqref{homo-ibvp} we expect the solution to decay   exponentially as $t\to\infty$. This can be corroborated in two ways: first, by numerically solving the problem directly; and second, by evaluating the first few modes of the general
``classical'' sine series solution representation
\begin{flalign}
\label{utm-series}
&u(x, t) 
=
\frac{1}{D}
\sum_{n=1}^\infty
e^{-(\frac{n\pi}{2D})^2 t}
\ \hats{u}_{\!\!\!0}\big(\tfrac{n\pi}{2D}\big)
 \sin\!\left(\tfrac{n\pi(x+D)}{2D}\right)&
\nonumber\\
&
+
\frac{1}{D}
\sum_{n=1}^\infty
e^{-(\frac{n\pi}{2D})^2 t}
\Big\{
\tfrac{n\pi}{2D} 
\Big[
\tilde g_0\big(\!\!\left(\tfrac{n\pi}{2D}\right)^2\!, t\big)
 -
e^{in\pi }\tilde h_0\big(\!\!\left(\tfrac{n\pi}{2D}\right)^2\!,t\big)
\Big]\!\Big\}&
\nonumber\\
&\hskip 4.7cm\times \sin\!\left(\tfrac{n\pi(x+D)}{2D}\right),&
\end{flalign}
where $\hats{u}_{\!\!\!0}$ denotes the sine transform of $u_0$ 
on $[-D, D]$ defined by
\begin{equation}\label{sine-trans-def}
\hats{u}_{\!\!\!0}(k)= \int_{-D}^D \sin\left(k(x+D)\right)u_0(x) dx,
\end{equation}
and the vanishing -- in this case -- terms involving
$\tilde{g}_0$ and $\tilde{h}_0$ are defined by Eq.~\eqref{gh0til-def}.
The outcome of this analysis is shown in Fig.~\ref{fig2}, which is
clearly showcasing the exponential decay of the solution of
this problem.

\begin{figure}[tbp]
\begin{center}
\includegraphics[width=6.5cm]{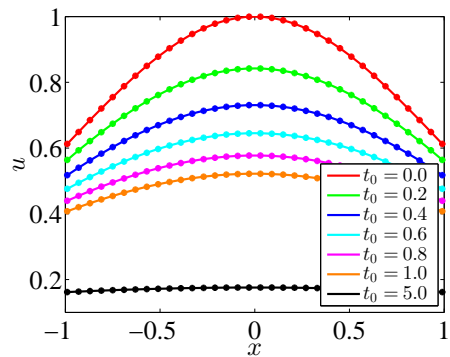}
\caption{
Spatial profiles of the solution $u$ of 
the IBVP with initial condition $u(x,0)=e^{-x^2/2}$ and 
boundary conditions given by Eq.~\eqref{kaz1} at various
times $t=t_{0}$. The asterisks denote the solution $u$
obtained by direct numerical computation,
while the solid lines the corresponding semi-analytical
approximation. 
In the latter, the decomposition~\eqref{kaz2} has been used
with the component $u^{(2)}$ given by 
the exact formula~\eqref{kaz3} and
the component $u^{(1)}$ 
obtained by using the first ten modes of the series 
representation~\eqref{utm-series}.
}
\label{fig22}
\end{center}
\end{figure}

Previously, in Fig.~\ref{case4}, we verified that the transparency
induced by consonant and compatible problems enables arbitrary power law
decays for IBVPs. Now, in Fig.~\ref{fig22}, through the accuracy
of the comparison of the numerical solution with the decomposition
(into problems~\eqref{homo-ibvp} and~\eqref{kummer-ibvp}) of compatible but non-consonant IBVPs,
we confirm that the solution of those problems is not genuinely self-similar,
yet it can be decomposed in a self-similar power law decay and
an exponential one associated with homogeneous boundary conditions.

\begin{figure}[tbp]
\centering
\hskip -1.5cm
\includegraphics[width=7cm]{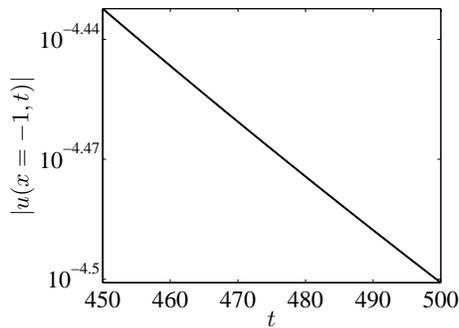}
\caption{
Semi-logarithmic plot of $|u(-1, t)|$, 
where $u(x, t)$ is the solution of the diffusion equation~\eqref{deq1}
on the interval $[-D, D]$ with initial condition $u_0(x) = x e^{-x^2/2}$,
boundary conditions $u(\pm D, t) = \pm \frac{D}{(2t+1)^{3/2}} e^{-\frac{D^2}{2(2t+1)}}$
and $D=200$. Since the boundary conditions are smaller than machine precision, 
the solution decays exponentially rather than algebraically.}
  \label{case4-large}
\end{figure}

As a final side remark and a note of caution
stemming from direct numerical observations,
we should point out that the case of boundary conditions that are both compatible and consonant with the initial condition can nevertheless
be sensitive to the choice of the size of the spatial domain $[-D, D]$. 
Indeed, although in this case one would normally expect to observe the solution decaying at the \textit{algebraic} rate imposed by the initial data and ``preserved'' by the boundary data (cf. Fig.~\ref{case4}), if the choice of $D$ is such that the size of the boundary data is smaller than machine precision then one is effectively solving a \textit{homogeneous} Dirichlet problem. In this scenario, one will therefore observe \textit{exponential}, as opposed to algebraic, decay.
 For example, revisiting the IBVP of Fig.~\ref{case4} but now with $D=200$ instead of $D=1$ makes the boundary conditions $u(\pm D, t) = \pm \frac{D}{(2t+1)^{3/2}} e^{-\frac{D^2}{2(2t+1)}}$ effectively zero. Hence, the decay rate observed for the solution is exponential instead of algebraic, as shown in Fig.~\ref{case4-large}.

\section{Conclusions and Future Work}

In the present work we have revisited the topic of self-similar
solutions of one of the most fundamental models, namely the linear
diffusion equation. Employing the 
MN-dynamics approach, we have justified  the existence of a broad class of self-similar
solutions and have developed a superposition
of the relevant eigenfunctions  upon performing a separation
of variables in the MN-frame. This has enabled us to identify solutions that decay exponentially 
as well as solutions that decay with power laws in the self-similar spatial variable.
It has also provided us with eigenfunctions that are bounded,
as well as others that are integrable under suitable power law decay
conditions (and even ones that are non-integrable).
Among the ones that are integrable (and hence physically
correspond to finite mass), the slowest
one identified has been the customary $t^{-1/2}$ decay law.
This is suggestive about the observability of this type of
decay in physical experiments. Nevertheless, we showcased explicit
examples featuring a different type of decay, arising from initial
data with vanishing projection on the corresponding eigenfunction.
In this way, we effectively showed that arbitrary power laws as
regards the temporal decay are
possible to realize in the context of the linear diffusion model.

On the other hand, we also revisited initial-boundary value problems (IBVPs) associated with the diffusion equation and identified different
possibilities. We proposed the notion of consonant boundary conditions,
for which a self-similar solution remains at all times transparent
to the presence of the boundaries and preserves its corresponding rate
of decay. This, in turn, illustrates that the associated IBVP can also feature
arbitrary temporal power law decay rates. On the other hand, compatible, yet
non-consonant IBVPs were formulated with different decay rates prompted
by their initial and boundary data. In this case, we established
that the general solution is not self-similar, but instead can be
decomposed in a self-similar part and a part with exponential decay
due to its satisfying homogeneous boundary conditions (this is the
typical, widely explored case considered in the context of
Fourier series).

These considerations pave a number of directions for future studies.
In the context of the diffusion model, it is well-known (see for
details Appendix A) that an analytical solution exists either
in the form of the \textit{unified transform method} of Fokas~\cite{f1997,fbook}, or in the form of an infinite
series representation. Such formulae already encompass a decomposition
of the solution into a component that stems from the initial condition
and a component that arises from the boundary contributions. In this context, it would
be especially interesting  to reconstruct
from these analytical expressions the parts of the solution that
may feature self-similarity, as well as to identify the ones that
do not, and to formulate more precise conditions under which we should
expect the solution to be self-similar. On the other hand, a far
more open-ended problem stems from the introduction of nonlinearity
in the system. Here, a principal question
concerns the potential persistence, as well as the
(appropriately adapted notion of) stability of the self-similar solutions
over the time evolution. An example in this direction is, for instance, given
in the piecewise linear model of~\cite{sr1998}.
Lastly, an especially important technical question even for the linear
diffusion case concerns the properties of the operator~\ref{deq9}
and potential decompositions of general solutions of this PDE on 
a basis of associated eigenfunctions.
Such questions are presently under
investigation and will be reported in future publications.

\section*{\uppercase{Appendix A: Analytical Solution
    of the Diffusion IBVP using the Unified Transform Method of
  Fokas}}

\begin{figure}[ht]
\begin{center}
\includegraphics[scale=.9]{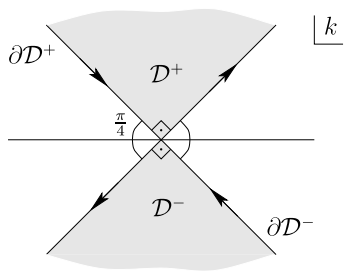}
\caption{
The regions $\mathcal D^+$ and $\mathcal D^-$ and the positively oriented boundaries $\partial \mathcal D^+$ and $\partial \mathcal D^-$.
}
\label{dpm}
\end{center}
\end{figure}

With the question of compatible, yet non-consonant
boundary data for IBVPs, it is useful to recall that the diffusion
equation \eqref{deq1} formulated on the  interval 
with any admissible combination of boundary conditions can be solved
analytically via the unified transform method (UTM) of Fokas~\cite{f1997, fbook}.
For example, in the case of Dirichlet boundary conditions $u(-D, t)=g_0(t)$ and $u(D, t)=h_0(t)$  
with initial condition $u(x, 0)=u_0(x)$, the UTM yields the solution formula 
\begin{align}\label{utm-sol}
&u(x, t)
=
\frac{1}{2\pi}\int_{k\in\mathbb R} e^{ikx-k^2t} \widehat u_0(k) dk
\\
& 
-
\frac{1}{2\pi}\int_{k\in\partial \mathcal D^+} 
\frac{e^{ik(x+D)-k^2 t}}{e^{2ikD}-e^{-2ikD}}
\Big[-2ik e^{-2ikD} \tilde g_0(k^2, t)
\nonumber\\
&\hskip 1.32cm
+2ik \tilde h_0(k^2,t)+e^{ikD}\widehat u_0(k)-e^{-ikD} \widehat u_0(-k)\Big]dk
\nonumber\\
& 
-
\frac{1}{2\pi}\int_{k\in\partial \mathcal D^-} 
\frac{e^{ik(x-D)-k^2 t}}{e^{2ikD}-e^{-2ikD}}
\Big[-2ik  \tilde g_0(k^2, t)
\nonumber\\
&\hskip 0.5cm
+2ik e^{2ikD} \tilde h_0(k^2,t)+ e^{-ikD}\widehat u_0(k)- e^{ikD}\widehat u_0(-k)\Big]dk,
\nonumber
\end{align}
where  the transforms $\widehat u_0, \tilde g_0, \tilde h_0$ of the prescribed initial and boundary data are defined by
\begin{subequations}
\begin{align}
\widehat u_0(k) &= \int_{-D}^D e^{-ikx} u_0(x) dx, 
\label{u0hat-def}
\\
\big(\tilde g_0, \tilde h_0\big) (k^2, t) &= \int_0^t e^{k^2 t'} \left(g_0, h_0\right)(t')dt',
\label{gh0til-def}
\end{align}
\end{subequations}
and the complex contours of integration $\partial \mathcal D^\pm$ are the positively oriented boundaries of the regions $\mathcal D^\pm$ shown in Fig.~\ref{dpm}. 

Furthermore, the UTM solution of any IBVP for the diffusion equation on the interval 
can always be reduced to the corresponding ``classical'' infinite series solution representation \eqref{utm-series}.

\acknowledgments We are grateful to C. Eugene Wayne for numerous enlightening discussions on the subject and useful comments on the manuscript.
P.G.K. gratefully acknowledges the support of NSF-DMS-1312856, 
NSF-PHY-1602994 and the ERC under FP7, Marie Curie Actions, People,
International Research Staff Exchange Scheme (IRSES-605096). 
The research of I.G.K. was partially supported by the US  National Science Foundation and by
the IAS-TU Munich  through a Hans Fischer Senior Fellowship.
D.M. and E.G.C. would like to thank the Department of Chemical and Biological Engineering at Princeton University for the kind hospitality.

\end{document}